\documentclass[a4paper,11pt]{article}
\usepackage{a4wide}
\usepackage{theorem}
\usepackage{amsmath}
\usepackage{array}
\usepackage{amssymb}
\usepackage{amsfonts}
\usepackage[french]{babel}
\usepackage{epsf}
\usepackage{epsfig}

\newtheorem{theo}{\indent Theorem\newline}[section]

{\theorembodyfont{\rmfamily}%
\newtheorem{rem}[theo]{\noindent Remark}}
{\theorembodyfont{\rmfamily}%
 \theoremstyle{break}%
}
\newtheorem{prop}[theo]{\indent Proposition\newline}
\newtheorem{lemma}[theo]{\indent Lemma\newline}
\newtheorem{cor}[theo]{\indent Corollary\newline}

\def\Z{{\Bbb{Z}}}

\def\R{{\Bbb{R}}}
\def\C{{\Bbb{C}}}

\newcommand{\Div}{\mathop{\rm div}\nolimits}
\newcommand{\Jac}{\mathop{\rm Jac}\nolimits}
\newcommand{\Pic}{\mathop{\rm Pic}\nolimits}
\newcommand{\conj}{\mathop{\rm conj}\nolimits}

\newlength{\indentation}%
\makeatletter
\newcommand\@makefntextsans[1]{%
    \parindent 0em%
    \noindent%
    \hb@xt@0em{\hss}%
    #1}
\def\footnotetextsans{%
     \@ifnextchar [\@xfootnotenextsans%
       {\@footnotetextsans}}
\def\@xfootnotenextsans[#1]{%
  \begingroup%
     \csname c@\@mpfn\endcsname #1\relax%
  \endgroup%
  \@footnotetextsans}
\long\def\@footnotetextsans#1{\insert\footins{%
    \reset@font\footnotesize%
    \interlinepenalty\interfootnotelinepenalty%
    \splittopskip\footnotesep%
    \splitmaxdepth \dp\strutbox \floatingpenalty \@MM%
    \hsize\columnwidth \@parboxrestore%
    \color@begingroup%
      \@makefntextsans{%
        \rule\z@\footnotesep\ignorespaces#1\@finalstrut\strutbox}
    \color@endgroup}}
\makeatother

\begin{document}

\cleardoublepage
\title{Deformation classes of real ruled manifolds}
\author{Jean-Yves Welschinger }
\date{}
\maketitle

\makeatletter\renewcommand{\@makefnmark}{}\makeatother
\footnotetextsans{Keywords :  Ruled manifold, real algebraic manifold, 
deformation.}
\footnotetextsans{AMS Classification : 14J10, 14M99, 14P25.}

{\bf Abstract :}

A complete description of the deformation classes of real ruled manifolds 
is given. In particular, we prove that once the complex deformation class
is fixed, the real deformation class is prescribed by the topology of the
real structure.

\section*{Introduction}

A {\it real algebraic manifold} $(X , c_X)$ is a smooth complex algebraic 
manifold
$X$ equipped with an antiholomorphic involution $c_X$. The {\it real part} of
$X$ is the fixed point set of $c_X$. One of the main problems in real 
algebraic geometry nowadays is to understand the deformation classes 
(see \S \ref{subsectionoverc} for the definition) of
real algebraic manifolds. One can think of this problem as a modern version
of a question of Hilbert in his $16^{th}$ problem concerning the topology of
smooth real quartics in the real projective $3$-space. 
Several works
have already been done to solve this problem : the cases of real curves,
real rational surfaces, real minimal ruled surfaces and real minimal surfaces 
of Kodaira dimension $0$
are known (see \cite{Nat2}, \cite{KhDg2}, \cite{Wels}, \cite{KhDg},
\cite{Cat} and \cite{DIK} for an extension of \cite{KhDg} to
finite group actions on $K_3$'s). The purpose of this paper is to extend the 
result of \cite{Wels}
to ruled manifolds of higher dimensions. Note that the complex deformation 
classes of ruled manifolds were studied in \cite{Schm}.

A {\it ruled manifold} is a smooth algebraic manifold equipped with a proper
holomorphic submersion on a smooth compact irreducible curve $B$, whose fibers
are projective spaces. In dimension two,
these are geometrically ruled surfaces. The aim of this paper is to
prove the following theorem :

\begin{theo}
Two real ruled manifolds are in the same real deformation class if and only if
they are in the same complex deformation class and they are diffeomorphic via 
an equivariant diffeomorphism. 
\end{theo}
Moreover, all the deformation classes of real ruled manifolds will be 
described (see \S \ref{subsectiondefr}). It is necessary here to fix the 
complex deformation class of the manifold since as was noticed by E. Brieskorn
(see \cite{Brie}, satz $3.1$), there exist complex ruled manifolds which are 
diffeomorphic
to each other but which are not deformation equivalent. However, once
the complex deformation class is fixed, the topology of the involution is 
enough to describe the real deformation classes of the real manifolds, which 
is the case in all the known examples nowadays. 

The paper is organized as follows : in the first section, we give basic facts
and preliminary results on ruled manifolds, real structures on these
manifolds and a notion of elementary transformations that can be performed
 on them.
The second section is devoted to the statements of the results and the
third to their proofs.

\section{Real ruled manifolds and elementary transformations}

\subsection{Ruled manifolds}
\label{subsubsectionruledmfd}

A smooth irreducible compact complex manifold $X$ of dimension $n$ is
said to be {\it ruled} if there exists a smooth irreducible compact complex
curve $B$ and a proper holomorphic submersion $p : X \to B$ whose fibers are 
isomorphic to the projective space $\C P^{n-1}$.
For example, let $E$ be a complex vector bundle of rank $n$ over the curve $B$
and $X = P(E)$ be the associated projective bundle. Then $X$ is a ruled 
manifold.
Note that when $n=2$, ruled manifolds are geometrically ruled surfaces.
For these surfaces, it is well known that the curve $B$ is
unique, so as the ruling $p$ except from $X = \C P^1 \times \C P^1$ (see
\cite{Beau}). The following lemma extends this result.

\begin{lemma}
\label{lemmabase}
Let $X$ be a ruled manifold of dimension $n \geq 3$. Then the ruling $p$ and 
the curve $B$ are unique.
\end{lemma}

The curve $B$ is called the {\it base} of $X$.\\

{\bf Proof :}

It follows from the fact that the only divisors of $X$ isomorphic to 
$\C P^{n-1}$ are the fibers of $p$. Indeed, $p$ would restrict otherwise
to a surjective morphism from this divisor onto $B$, and the generic fibers of
this morphism would give smooth disjoint complex hypersurfaces of
$\C P^{n-1}$. Such hypersurfaces do not exist in dimension $n-1 \geq 2$.
$\square$.\\

The following proposition is mentioned in \cite{Schm}, p. $214$.

\begin{prop}
\label{proppe}
Let $X$ be a ruled manifold of dimension $n$ over $B$. Then there exists
a complex vector bundle $E$ of rank $n$ over $B$ such that $X$ is 
isomorphic to the projective bundle $P(E)$. Moreover, the projective bundle
$P(E')$ is isomorphic to $P(E)$ if and only if $E' = E \otimes L$ for
$L \in \Pic (B)$. $\square$
\end{prop}

\begin{cor}
Ruled manifolds are all projective algebraic. $\square$
\end{cor}

\begin{rem}
Let $L \in \Pic (B)$ and $E$ be a complex vector bundle of rank $n$ over $B$.
Then $\deg (E \otimes L) = \deg (E) + n \deg (L)$, where $\deg (E)$ stands
for the degree of $E$.
\end{rem}

Let $X = P(E)$ be a ruled manifold of dimension $n$ over $B$. We define the
{\it degree of $X$} to be $\deg (E)$ reduced modulo $n$. It will be denoted 
by $\deg (X) \in \Z / n \Z$.
Let $L$ be a complex line bundle over $B$ and $L_0$ be the trivial line 
bundle. The section $P(L)$ (resp. $P(L_0)$) of the ruled surface 
$P (L \oplus L_0)$ defines a divisor on this surface denoted by $D_L$
(resp. $D_{L_0}$).

\begin{lemma}
\label{lemmanorm}
1. Let $p$ be the ruling $P (L \oplus L_0) \to B$ and ${\cal O} (D_L)$,
${\cal O} (D_{L_0})$ denote the invertible sheaves associated to
the divisors $D_L$ and $D_{L_0}$. Then
$${\cal O} (D_L) = {\cal O} (D_{L_0}) \otimes p^* (L^*).$$

2. Let $F$ be a complex vector bundle over $B$, $X$ be the ruled manifold
$P (L \oplus L_0 \oplus F)$ and $N$ be the normal bundle of $P (L \oplus L_0)$
in $X$. Then :
$$N =  p^* (F) \otimes {\cal O} (D_L).$$
\end{lemma}

{\bf Proof :}

Let $D =
\sum_{i=1}^k n_i p_i $ be a divisor associated to $L$, where $p_i 
\in B$ and $n_i \in \Z$ for $i \in \{1, \dots , k\}$. Denote by $U_0 = B
\setminus \{ p_i \, | \, 1 \leq i \leq k \}$ and for every $i \in \{1,
\dots , k\}$, choose some holomorphic chart $(U_{p_i} , \phi_{p_i})$
such that $U_{p_i} \cap U_{p_j} = \emptyset$ if $i \neq j$ and 
$\phi_{p_i} : U_{p_i} \to \Delta = \{ z \in
\C \, | \, |z| < 1 \}$ is a biholomorphism satisfying
$\phi_{p_i} (p_i) = 0$.
For every $i \in \{1, \dots , k\}$, denote by $\psi_i$ the morphism :

$$\begin{array}{rcl}
(U_{p_i} \setminus p_i) \times \C P^1 & \to & U_0 \times \C P^1\\
(x,(z_1 : z_0)) & \mapsto & (x , (\phi_{p_i} (x)^{-n_i} z_1 : z_0)).
\end{array}$$
The morphisms $\psi_i$ allow to glue together the trivializations
$U_{p_i} \times \C P^1$, $i \in \{0, \dots , k\}$, in order to define the
ruled surface $P (L \oplus L_0)$. 

Let $f : p^{-1} (U_0) \to \C$, $(x,(z_1 : 1)) \mapsto z_1$. Then $f$
extends to a meromorphic function on $P (L \oplus L_0)$ such that
$f^{-1} (0) = D_L + \sum_{n_i \geq 0} n_i p^{-1} (p_i)$ and
$f^{-1} (\infty) =  D_{L_0} + \sum_{n_i \leq 0} n_i p^{-1} (p_i)$.
Hence $\Div (f) = D_L -  D_{L_0} + p^{-1} (D)$, so that
${\cal O} (D_L) = {\cal O} (D_{L_0}) \otimes p^* (L^*)$, which proves
the first part of the lemma.

To prove the second part of the lemma, take a refinement of the covering
$(U_i)$ such that the bundle $F$ is trivial over any element of this covering.
The manifold $X$ is then defined as the gluing of charts ${\cal U}_i^0$
and ${\cal U}_i^1$ isomorphic to $U_i \times \C^{n-1}$ with gluing maps :
${\cal U}_i^0 \cap {\cal U}_i^1 \to {\cal U}_i^1 \cap {\cal U}_i^0$,
$(x,z_1 ,f) \mapsto (x,\frac{1}{z_1} ,\frac{1}{z_1} f)$ where $f$ is a local
trivialization of $F$ over $U_i$, and ${\cal U}_i^0 \cap {\cal U}_j^0 \to 
{\cal U}_j^0 \cap {\cal U}_i^0$, $(x,z_0 ,f) \mapsto (x, l_{ij}^{-1} z_0,
l_{ij}^{-1} g_{ij} (f))$ where $l_{ij}$ and $g_{ij}$ are the changes of 
trivialization of $L$ and $F$ respectively. We deduce that the normal bundle
$N$ of $P (L \oplus L_0)$ in $P (L \oplus L_0 \oplus F)$ is defined as
the gluing of the trivializations :
$\big( {\cal U}_i^0 \cap {\cal U}_i^1 \cap P (L \oplus L_0) \big) 
\times \C^{n-2} \to \big( {\cal U}_i^1 \cap {\cal U}_i^0 \cap P (L \oplus L_0)
 \big) \times \C^{n-2}$,
$((x,z_1) ,\nu) \mapsto ((x,\frac{1}{z_1}) ,\frac{1}{z_1} \nu)$ and
$\big({\cal U}_i^0 \cap {\cal U}_j^0 \cap P (L \oplus L_0) \big) \times 
\C^{n-2} \to 
\big( {\cal U}_j^0 \cap {\cal U}_i^0 \cap P (L \oplus L_0) \big)\times 
\C^{n-2}$, 
$((x,z_0) ,\nu) \mapsto ((x, l_{ij}^{-1} z_0),
l_{ij}^{-1} g_{ij} (\nu))$.

Hence $N = p^* (F) \otimes {\cal O} (D_{L_0}) \otimes p^* (L^*) =
p^* (F) \otimes {\cal O} (D_L)$. $\square$

\begin{prop}
\label{propaut}
Let $B$ be a smooth irreducible compact complex curve and $L$ be a complex
line bundle over $B$ such that $L \neq L^*$ if $L$ is non-trivial. Let
$E = (L \oplus L_0)^k$ and $X = P(E)$. Then every automorphism of $X$ fibered 
over the identity of $B$ which leaves the $k$ ruled surfaces $P(L \oplus L_0)$
invariant lifts to an automorphism of $E$ fibered over the 
identity of $B$.
\end{prop}

{\bf Proof :}

Let $\phi$ be such an automorphism of $X$. From proposition $2.1$ of
\cite{Wels}, we know that the restriction of $\phi$ to the $j^{th}$ ruled
surface $P(L \oplus L_0)$ lifts to an automorphism $\psi_j$ of the rank
two vector bundle $L \oplus L_0$. Let $\Psi$ be the automorphism
$(\psi_1 , \dots , \psi_k)$ of $E$. This automorphism induces an automorphism 
$\psi$ of $X$ such that $\psi^{-1} \circ \phi$ is the identity once
restricted to each ruled surface $P(L \oplus L_0)$. But such an automorphism
lifts to a diagonal automorphism of $E$ of the form $(\lambda_1 Id, \dots ,
\lambda_k Id)$, where $\lambda_j \in \C^*$.
Hence the result. $\square$

\subsection{Real structures on ruled manifolds}
\label{subsubsectionrealruledmfd}

A {\it real structure} on the ruled manifold $X$ is an antiholomorphic 
involution $c_X : X \to X$. The fixed point set of $c_X$ is called the 
{\it real part} of $X$ and is denoted by $\R X$.

\begin{lemma}
Let $p :X \to B$ be a ruled manifold of dimension $n >2$ and $c_X$ be a real
structure on $X$. Then there exists a real structure $c_B$ on $B$ such
that $p \circ c_X = c_B \circ p$. Moreover, this real structure $c_B$ is
unique.
\end{lemma}

The real structure $c_X$ will be said to be {\it fibered} over $c_B$.\\

{\bf Proof :}

From the proof of lemma \ref{lemmabase}, we know that the only divisors of $X$ 
isomorphic to $\C P^{n-1}$ are the fibers of $p$. So $c_X$ preserves these 
fibers and hence induces a diffeomorphism $c_B$ on the base. This 
diffeomorphism
is antiholomorphic and is an involution. $\square$\\

We deduce from this lemma that the connected components of $\R X$ are
$\R P^{n-1}$-bundles over the circle. For odd $n$,
such a bundle is unique whereas for even $n$ there are two such bundles,
one which is orientable and the other one which is not. We define the
 {\it topological 
type} of a real ruled manifold $(X , c_X)$ of even complex dimension $n$ to
be the quintuple of integers $(t,k,g, \mu, \epsilon)$ where $t$ is the
number of orientable components of $\R X$, $k$ is the number of non-orientable 
components of $\R X$ and $(g, \mu, \epsilon)$ is the topological type of
the real curve $(B , c_B)$, that is the genus of $B$, the number of connected 
components of $\R B$ and the dividing or non-dividing type of $(B , c_B)$. 
This definition extends the one given in \cite{Wels} for $n=2$.

Let us present now an important example of real ruled manifold.
Let $(B , c_B)$ be a real algebraic curve and $L$ be a complex line bundle 
over $B$ such that $c_B^* (L) = L^*$ where $c_B^*$ is the real structure on
$\Pic (B)$ induced by $c_B$ (see \cite{Wels}, \S $1.1$). Let $D$ be a divisor
associated to $L$ and $f_D$ be a meromorphic
function on $B$ such that $\Div (f_D) = D + c_B (D)$ and $f_D
= \overline{f_D \circ c_B}$ (it always exists, see \cite{Wels}, lemma 1.3).
Note that the sign of $f_D$ is constant on every
component of $\R B$. The following proposition is analogous to proposition
1.6 of \cite{Wels} :

\begin{prop}
\label{propcfd}
Associated to every such couple 
$(D, f_D)$ on $(B , c_B)$, there exists a real structure $c_{f_D}$ on 
$X = P ((L \oplus L_0)^\frac{n}{2})$ fibered over $c_B$,  
whose real part is orientable and maps surjectively onto the components 
of $\R B$ on which $f_D$ is non-negative.
\end{prop}

\begin{rem}
When there will not be any ambiguity on the choice
of the couple $(D, f_D)$, we will denote by $c_X^+$ (resp. $c_X^-$) the
real structure $c_{f_D}$  (resp. $c_{-f_D}$).
\end{rem}

{\bf Proof :}

Denote $D =
\sum_{i=1}^k n_i p_i $ where $p_i 
\in B$ and $n_i \in \Z$ for $i \in \{1, \dots , k\}$. We can assume that 
the set $\{ p_i
\, | \, 1 \leq i \leq k \}$ is invariant under $c_B$. Let $U_0 = B
\setminus \{ p_i \, | \, 1 \leq i \leq k \}$ and $(U_{p_i} ,
\phi_{p_i})$, $i \in \{1,
\dots , k\}$, be an atlas
compatible with the divisor $D$ and the group $<c_B>$ (see \cite{Wels}, 
page 3). 

The morphisms :

$$\begin{array}{rcl}
(U_{p_i} \setminus p_i) \times \C P^{n-1} & \to & U_0 \times \C P^{n-1} \\
(x,(z_1^j : z_0^j)) & \mapsto & (x , (\phi_{p_i} (x)^{-n_i} z_1^j : z_0^j))
\end{array}$$
($i \in \{1, \dots , k\}$, $j \in \{1, \dots , \frac{n}{2} \}$) allow to glue 
together the trivializations
$U_{p_i} \times \C P^{n-1}$, $i \in \{0, \dots , k\}$, in order to define the
ruled manifold $X$. 

Now, the maps 
$$\begin{array}{rcl}
U_0  \times \C P^{n-1}& \to & U_0 \times \C P^{n-1} \\
(x,(z_1^j : z_0^j)) & \mapsto & (c_B (x) , (\overline{z_0^j} : f_D \circ c_B
(x) \overline{z_1^j})), 
\end{array}$$
and for every $i \in \{1, \dots , k\}$,
$$\begin{array}{rcl}
U_{p_i}  \times \C P^{n-1}& \to &  U_{c_B(p_i)} \times \C P^{n-1}\\
(x,(z_1^j : z_0^j)) & \mapsto & (c_B (x) , (\overline{z_0^j} : f_D \circ c_B
(x)\overline{\phi_{p_i} 
(x)}^{-n_{c_B(p_i)}-n_{p_i}} \overline{z_1^j}),
\end{array}$$
where $j \in \{1, \dots , \frac{n}{2} \}$, glue together to form an 
antiholomorphic map $c_{f_D}$ on $X$. This map
lifts $c_B$ and is an involution. 

Now, the fixed point set of $c_{f_D}$ in $U_0 \times \C
P^{n-1}$ is :
$$\{ (x, (\overline{\theta_j} : \sqrt{f_D (x)}\theta_j)) \in U_0 \times \C
P^{n-1} \, | \, x \in \R B, f_D (x) \geq 0 \text{ and } \theta_j \in \C 
\text{ for }
j \in \{1, \dots , \frac{n}{2} \}  \}.$$
The connected components of this fixed point set are orientable 
$\R P^{n-1}$-bundle over a circle or an interval depending on whether 
the corresponding component of $\R B$ is
completely included in $U_0$ or not. Similarly, the fixed point set of
$c_{f_D}$ in $U_{p_i} \times \C P^{n-1}$ is :
$$\{ (x, (\overline{\theta^i_j} : \sqrt{f_D (x) \times x_i^{-2n_i}}
\theta^i_j)) \in U_{p_i} \times \C
P^{n-1} \, | \, x \in \R B, f_D (x) \geq 0 \text{ and } \theta^i_j \in \C ,
j \in \{1, \dots , \frac{n}{2} \}   \},$$ 
where $x_i = \phi_{p_i} (x)$. This fixed point set is a cylinder if
$p_i \in \R B$ and is empty otherwise.

The gluing maps between these cylinders are given by $\theta_j =
\sqrt{-1} \theta^i_j$ if $x_i = \phi_{p_i} (x) < 0 $ and by $\theta_j =
\theta^i_j$ if $x_i = \phi_{p_i} (x) > 0 $. Since these two maps
preserve the orientation of $\R P^{n-1}$, the results of these gluings 
are always orientable.
Thus, the real part of $(X, c_{f_D})$ consists only of orientable components
and these components stand exactly over the components of $\R B$ 
on which $f_D \geq 0$. $\square$

\begin{prop}
\label{propconj}
Let $(B , c_B)$ be a real algebraic curve with non-empty real part and 
$L \in \Pic (B)$ be such that $c_B^* (L) = L^*$ and $L \neq L^*$ if $L$
is non-trivial. Let $n$ be an even integer, 
$X = P ((L \oplus L_0)^\frac{n}{2})$  and $c_X$ be
a real structure on $X$ which leaves the $\frac{n}{2}$ ruled surfaces
$P (L \oplus L_0) \subset X$ invariant. Then $c_X$ is conjugated to one of 
the two real structures $c_X^\pm$.
\end{prop}

{\bf Proof :}

If $L$ is trivial, then $X = B \times \C P^{n-1}$ and $c_X^\pm$ are the
two real structures $c_B \times \conj$ and $c_B \times c_0$, $c_0$ being
the real structure of $\C P^{n-1}$ with empty real part. Now the automorphisms
of $X$ fibered over the identity of $B$ are just the automorphisms of 
$\C P^{n-1}$, so every real structure
of $X$ fibered over $c_B$ is the product of $c_B$ with a real structure of 
$\C P^{n-1}$. Hence the result follows from
the well known fact that the standard complex conjugation $\conj$ and $c_0$
are the only real structures of $\C P^{n-1}$ up to conjugation. Let us now
assume that $L$ is non-trivial. The real structure $c_X$ can be written
$c_X^+ \circ \phi$, where $\phi$ is an automorphism of $X$ fibered over the
identity of $B$ and which leaves the $\frac{n}{2}$ ruled surfaces
$P (L \oplus L_0) \subset X$ invariant. From proposition \ref{propaut}, there
exists an automorphism $\Phi$ of the vector bundle $(L \oplus L_0)^\frac{n}{2}$
which lifts $\phi$. Since $L$ is non-trivial, $H^0 (B ; L) = H^0 (B ; L^*)=0$,
so that $\Phi$ is diagonal of the form $(a_1 , b_1 , a_2 , b_2, \dots ,
a_{\frac{n}{2}} , b_{\frac{n}{2}})$ with $a_j , b_j \in \C^*$. Since
$c_X^2$ is the identity, $a_1 \overline{b_1} = a_2 \overline{b_2} = \dots =
a_{\frac{n}{2}} \overline{b_{\frac{n}{2}}} \in \R^*$. Thus, dividing $\Phi$
by a real constant if necessary, we can assume that $a_1 \overline{b_1} = 
a_2 \overline{b_2} = \dots = a_{\frac{n}{2}} \overline{b_{\frac{n}{2}}} = 
\pm 1$. Moreover, we can assume that all the coefficients $a_j$ are equal to
one, replacing $c_X$ by its conjugated with the automorphism
$(a_1^{-1}, 1, a_2^{-1}, 1, \dots , a_{\frac{n}{2}}^{-1}, 1)$ otherwise.
Then, either all the coefficients $b_j$ are equal to $1$, or they are all 
equal to $-1$. In the first case, $c_X$ is conjugated to $c_X^+$, in the second
case, it is conjugated to $c_X^-$. $\square$

\subsection{Elementary transformations}

Let $X$ be a ruled manifold of dimension $n$ over the curve $B$. Let $x \in B$,
$X_x = p^{-1} (x)$ and $H_x , K_x \subset X_x$ be two disjoint projective
subspaces of $X_x$ of dimensions $k$ and $n-2-k$ respectively, where
$0 \leq k \leq n-2$. The blowing up $Y$ of $X$ along $H_x$ creates an
exceptional divisor $E_x$ isomorphic to $P (N_x \oplus L_0)$ where $N_x$ is
the normal bundle of $H_x$ in $X_x$. The strict transform $\widetilde{X}_x$
of $X_x$ in $Y$ intersects $E_x$ in the submanifold $P (N_x) \subset 
P (N_x \oplus L_0)$. Moreover, $\widetilde{X}_x$ is a ruled manifold over
$K_x$. Indeed, for every $y \in K_x$ consider the projective subspace of 
dimension $k+1$ of $X_x$ containing $H_x$ and passing through $y$. These
projective subspaces form a singular pencil parametrized by $K_x$. This
pencil lifts on $Y$ to a ruled manifold of base $K_x$ and fiber $\C P^k$,
which coincide with $\widetilde{X}_x$.
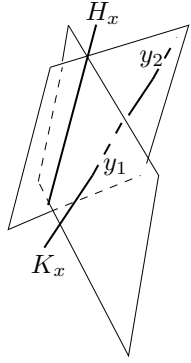
$$\vcenter{\hbox{\begin{picture}(0,0)%
\epsfig{file=rmfd1.pstex}%
\end{picture}%
\setlength{\unitlength}{3315sp}%
\begingroup\makeatletter\ifx\SetFigFont\undefined%
\gdef\SetFigFont#1#2#3#4#5{%
  \reset@font\fontsize{#1}{#2pt}%
  \fontfamily{#3}\fontseries{#4}\fontshape{#5}%
  \selectfont}%
\fi\endgroup%
\begin{picture}(1374,2712)(2689,-5923)
\put(3286,-3391){\makebox(0,0)[lb]{\smash{\SetFigFont{10}{12.0}{\rmdefault}{\mddefault}{\updefault}$H_x$}}}
\put(2881,-5281){\makebox(0,0)[lb]{\smash{\SetFigFont{10}{12.0}{\rmdefault}{\mddefault}{\updefault}$K_x$}}}
\put(3421,-4516){\makebox(0,0)[lb]{\smash{\SetFigFont{10}{12.0}{\rmdefault}{\mddefault}{\updefault}$y_1$}}}
\put(3691,-3706){\makebox(0,0)[lb]{\smash{\SetFigFont{10}{12.0}{\rmdefault}{\mddefault}{\updefault}$y_2$}}}
\end{picture}
}}$$

The composition of the blowing up of $X$ along $H_x$ and the blowing 
down of $\widetilde{X}_x$ on $K_x$ is called the 
{\it elementary transformation} of $X$ along $H_x$. For example if $n=2$,
then $H_x$ and $K_x$ are points in $X_x$ and the elementary transformation 
of $X$ along
 $H_x$ is the blowing up of the point $H_x$ composed with the blowing 
down of the strict transform of the fiber $X_x$. This notion of
elementary transformation thus extends the one used in \cite{Wels} for ruled 
surfaces.

\begin{lemma}
Let $\Delta = \{ z \in \C \, | \, |z| < 1 \}$, the elementary transformation of
$\Delta \times P (\C^{k+1} \times \C^{n - k-1})$ along $\{ 0 \} \times
P (\C^{k+1})$ is given by the following quadratic transform :
$$\begin{array}{rcl}
\Delta \times P (\C^{k+1} \times \C^{n - k-1})& \to &\Delta \times 
P (\C^{k+1} \times \C^{n - k-1}) \\
(x,(y_0 : \dots : y_k : z_0 :\dots : z_{n-2-k})) & \mapsto & 
(x , (t_0 : \dots : t_k : w_0 :\dots : w_{n-2-k}))
\end{array}$$
with $t_i = y_i$ for $0 \leq i \leq k$ and $xw_j = z_j$ for
$0 \leq j \leq n-k-2$.
\end{lemma}

{\bf Proof :}

It suffices to notice that the blowing up of $\Delta \times P (\C^{k+1} \times 
\C^{n - k-1})$ along $\{ 0 \} \times
P (\C^{k+1})$ embeds into $\Delta \times \C P^{n-1} \times \C P^{n-1}$. The
image of this embedding is :
$$Y = \big\{ (x,(y_0 : \dots : y_k : z_0 :\dots : z_{n-2-k}), (t_0 : \dots : 
t_k : w_0 :\dots : w_{n-2-k})) \in \Delta \times \C P^{n-1} \times \C P^{n-1} 
\, | \, $$
$$(y_0 : \dots : y_k) = (t_0 : \dots : t_k) \text{ or } y_i = t_i =0
\text{ and } (xy_i : z_0 :\dots : z_{n-2-k}) = (t_i : w_0 :\dots : w_{n-2-k})$$
$$\text{ for } 0 \leq i \leq k.\big\} $$
(The equations $(xy_i : z_0 :\dots : z_{n-2-k}) = (t_i : w_0 :\dots : 
w_{n-2-k})$ have to be ignored when they are not well defined.)
The projection of $Y$ onto the first factor $\Delta \times \C P^{n-1}_{(y,z)}$
(resp. onto the second factor $\Delta \times \C P^{n-1}_{(t,w)}$) is the
blowing up of $\Delta \times P (\C^{k+1} \times 
\C^{n - k-1})$ along $\{ 0 \} \times P (\C^{k+1})$ (resp. along
$\{ 0 \} \times P (\C^{n-k-1})$), hence the result. $\square$

\begin{cor}
\label{corelem}
Let $F,G$ be complex vector bundles over a curve $B$. Let $x \in B$ and
$X=P (F \oplus G)$. The ruled manifold obtained from $X$ after an
elementary transformation along $P (F_x)$ is the manifold $P (F(x) \oplus G)$
where $F(x) = F \otimes {\cal O}(x)$. $\square$
\end{cor}
The next corollary will be fondamental in what follows. It allows to break
into pieces or to glue together some elementary transformations.

\begin{cor}
\label{corcont}
Let $F,G,H$ be complex vector bundles over a curve $B$, $x \in B$ and
$X=P (F \oplus G \oplus H)$. The ruled manifold obtained from $X$ after an
elementary transformation along $P (F_x)$ followed by an
elementary transformation along $P (G_x)$ is the ruled manifold obtained from 
$X$ after an elementary transformation along $P (F_x \oplus G_x)$. $\square$
\end{cor}

\section{Statements of the results}

\subsection{Deformation over $\C$ of ruled manifolds}
\label{subsectionoverc}

Let $\Delta \subset \C$ be the Poincar\'e's disk equipped with the complex
conjugation $\conj$. Remember that a {\it deformation} of complex manifolds 
of dimension $n$ is a proper
holomorphic submersion $\pi : Y \to \Delta$ where $Y$ is an
analytic manifold of dimension $n+1$. If $Y$ is real and $\pi$ satisfies 
$\pi \circ c_Y =
conj \circ \pi$, then the deformation is said to be {\it real}. When 
$t \in ]-1 , 1 [ \in \Delta$, the fibers $Y_t =
\pi^{-1} (t)$ are invariant under $c_Y$ and are then compact real
analytic manifolds. Two complex (resp. real) analytic manifolds $X'$ and 
$X''$ are said to be 
{\it in the same deformation class} or {\it deformation equivalent} 
if there exists a chain $X'=X_0,
\dots , X_k=X''$ of compact complex (resp. real)
analytic manifolds such that for every $i \in \{ 0 , \dots , k-1 \}$, the
manifolds $X_i$ and $X_{i+1}$ are isomorphic to some complex (resp. real)
fibers of a complex (resp. real) deformation.

The following result can be found in \cite{Schm}.

\begin{theo}
Two complex ruled manifolds are in the same deformation class if and only if
they have same degree and bases of same genus. $\square$
\end{theo}
(see \S \ref{subsubsectionruledmfd} for the definition of the degree.)

Remember that ruled manifolds with opposite degrees and bases of same genus
are diffeormorphic even though they are not deformation equivalent, see
\cite{Brie}.

\subsection{Deformation over $\R$ of ruled manifolds}
\label{subsectiondefr}

\begin{theo}
\label{theoodd}
Two real ruled manifolds of odd dimension are in the same deformation class
if and only if they have same degree and bases of same
topological type. Similarly, 
two real ruled manifolds of even dimension are in the same deformation class
if and only if they have same degree, same
topological type and homeomorphic quotients.
\end{theo}

\begin{rem}
\label{remtheo}
For real ruled manifolds of even dimension, as soon as the real part of the 
base is non-empty, the condition on the quotients can be removed. However,
when the real part of the base is empty, there are two different deformation
classes of real ruled manifolds with same degree and same
topological type, see proposition \ref{propempty}. Note that when $n=2$, this
theorem \ref{theoodd} has already been obtained in \cite{Wels}, theorem $3.7$.
\end{rem}

Note that when $n$ is odd, the ruled manifolds 
$P(L(x + c_B (x))^d \oplus L_0^{n-d})$, where
$d \in \{0 , \dots , n-1 \}$, $x \in B \setminus \R B$ and $(B , c_B)$ is of
any topological type, are all real. Hence there do exist odd-dimensional
real ruled manifolds of any degree and any topological type (realized by a real
algebraic curve).
This together with theorem \ref{theoodd}
completely describes the deformation classes of real ruled manifolds of odd 
dimension.
A quintuple $(t,k,g, \mu ,\epsilon)$ of integers is called {\it allowable} 
if $t,k
\geq 0$, $t+k \leq \mu$ and $(g, \mu ,\epsilon)$ is the topological type
of a real curve. Obviously, the topological types of real ruled manifolds
(see \S \ref{subsubsectionrealruledmfd} for the definition) of even dimensions
are allowable.

\begin{prop}
\label{propeven}
Any allowable quintuple $(t,k,g, \mu ,\epsilon)$ is realized as the 
topological type of a real ruled manifold of any even dimension and any 
degree $d$ satisfying $d=k \mod (2)$. There do not exist such real ruled 
manifold of degree $d \neq k \mod (2)$.
\end{prop}

Note that in dimension two, this proposition has already been obtained in 
\cite{Wels}, proposition 3.4. Together with theorem \ref{theoodd} and
remark \ref{remtheo}, it completely describes the deformation classes of real 
ruled manifolds of even dimension. \\

{\bf Proof :}

Let $(t,k,g, \mu ,\epsilon)$ be an  allowable quintuple. There exists a smooth
compact connected real algebraic 
curve $(B,c_B)$ whose topological type is $(g, \mu ,
\epsilon)$ (see \cite{Nat2} for instance). If $\mu = 0$, the ruled manifold 
$(B \times \C P^{n-1} , c_B \times \conj)$, where $\conj$ is the standard real
 structure on $\C P^{n-1}$, is of topological type $(0,0,g,0,0)$. 
If $\mu \neq 0$, choose a partition ${\cal P}$ of $\R B$ in two elements
 such that one
of them contains $t+k$ components of $\R B$ and the other one $\mu - t
-k$. It follows from \cite{Wels}, lemma $3.2$, that there exists a
line bundle $L$ over $B$ such that $c_B^* (L) = L^*$ and the partition
associated to $L$ is ${\cal P}$ (see \cite{Wels}, \S $3.1$ for a definition). 
Thus, it follows from proposition \ref{propcfd}
that there exists a real structure $c_X^+$ on the ruled manifold 
$X=P((L \oplus L_0)^{\frac{n}{2}})$ such that the real part of $X$ consists
of $t+k$ orientable $\R P^{n-1}$ bundles. Choose $k$ of these orientable
bundles and make an elementary
transformation on a point of each of them. The result is a real ruled 
manifold of degree $k$ and topological type $(t,k,g, \mu ,\epsilon)$.
To get any other ruled manifold of degree congruent to $k$ modulo $2$
and same topological type, it suffices to perform the suitable number of
 elementary transformations 
along complex conjugated points. The fact that the condition $d=k \mod (2)$
is necessary will follow from the proof of theorem \ref{theoodd}. Indeed,
we will see that every real ruled manifold is deformation equivalent to a 
manifold obtained from a degree $0$ real ruled manifold with orientable real 
part after a finite number of elementary transformations performed on real
or complex conjugated points. Thus the degree modulo two of $X$ is encoded
by the topology of the real part, which finishes the proof. $\square$

\section{Proof of theorem \ref{theoodd}}

Since the case of real ruled manifolds having bases with empty real parts
requires special attention, the proof of theorem \ref{theoodd} in this
case is postponed to \S \ref{parempty}.

\subsection{When the base has non-empty real part}
\label{parnonempty}

The next two propositions will allow us to reduce the study of real ruled
manifolds to the study of some particular ones. Note that even if this 
paragraph is devoted to real ruled manifolds having bases with non-empty 
real parts, the assumption $\R B \neq \emptyset$ won't be made in
proposition \ref{propone} and lemma \ref{lemmadiv}.

\begin{prop}
\label{propone}
Let $(X , c_X)$ be a real ruled manifold over $(B , c_B)$. Then there
exists complex line bundles $L_1, \dots, L_k$ over $B$, complex vector bundles
$F_1, \dots, F_l$ of rank two over $B$ and a real structure $c_Y$ on
$Y=P(L_1 \oplus \dots \oplus L_k \oplus F_1 \oplus \dots \oplus F_l)$
such that :

- The sections $P(L_i)$, $i \in \{ 1, \dots , k\}$, and the ruled surfaces
$P(F_j)$, $j \in \{ 1, \dots , l\}$, of $(Y,c_Y)$ are all real.

- The real ruled surfaces $P(F_j) \subset Y$ , $j \in \{ 1, \dots , l\}$,
do not have any real holomorphic section.

- The real ruled manifolds $(X , c_X)$ and $(Y,c_Y)$ are in the same 
deformation class.
\end{prop}

The proof of this proposition is similar to the proof of proposition 3.8 in
\cite{Wels}.\\

{\bf Proof :}

Let $E$ be a complex vector bundle of rank $n$ over $B$ such that
$X=P(E)$ (it exists from proposition \ref{proppe}). If $X$ has a
real holomorphic section, denote by $M$ the sub-line bundle of $E$
associated to such a section. Otherwise, let $S$ be a holomorphic section
of $X$. Then $c_X (S) \neq S$, 
so these two sections intersect in a finite
number of points over the points $x_1, \dots , x_k$ of $B$ say. Let
$R \subset X$ be the ruled surface generated by $S$ and $c_X (S)$, that is
the closure of the surface whose fiber over $x \in B \setminus \{x_1, \dots , 
x_k \}$ is the line joining $S (x)$ to $c_X (S) (x)$ in $X_x$. By 
construction, $R$ is
a real ruled surface of $X$. Denote in this case by $M$ the rank two 
sub-bundle of $E$ associated to $R$. Now, let $N$ be the quotient bundle
$E/M$ so that $E$ is an extension of $N$ by $M$. Let $\mu \in H^1
(B , {\cal H}om (N,M))$ be the extension class of this bundle and let
$\mu^1$ be a $1$-cocycle with coefficients in the sheaf ${\cal H}om (N,M)$, 
defined on a covering ${\cal U}=(U_i)_{i \in I}$ of
$B$, realizing the cohomology class $\mu \in H^1
(B ,  {\cal H}om (N,M))$. The bundle $E$ is then obtained as the gluing of the
bundles $(M \oplus N)|_{U_i}$ by the gluing maps :

\begin{eqnarray*}
(M \oplus N)|_{U_i \cap U_j} & \to & (M \oplus N)|_{U_j \cap U_i} \\
(m,n) & \mapsto &  \left[ \begin{array}{cc}
1&\mu_{ij}\\
0&1
\end{array} \right]
\left( \begin{array}{c}
m\\
n
\end{array} \right)
=(m+\mu_{ij} n , n).
\end{eqnarray*}
We can assume that for every open set $U_i$ of ${\cal U}$, there
exists $\overline{\i} \in I$ such that $U_{\overline{\i}} = c_B (U_i)$
(add these open sets to ${\cal U}$ if not). We can also assume that
there exists $J \subset I$ such that the open sets $(U_i)_{i \in J}$
cover $B$ and such that the real structure $c_X : X|_{U_i} \to
X|_{U_{\overline{\i}}}$ lifts to an antiholomorphic map $E|_{U_i} \to
E|_{U_{\overline{\i}}}$ (take a refinement of ${\cal U}$ if
not). Since by hypothesis the section or ruled surface of $X$ associated to 
$M$ is
real, these antiholomorphic maps are of the form :
\begin{eqnarray*}
(M \oplus N)|_{U_i} & \to & (M \oplus N)|_{U_{\overline{\i}}} \\
(x,(m,n)) & \mapsto &  (c_B (x) , \left[ \begin{array}{cc}
a_i&b_i\\
0&d_i
\end{array} \right]
\left( \begin{array}{c}
m\\
n
\end{array} \right) ),
\end{eqnarray*}
where $a_i$ (resp. $b_i$, resp. $d_i$) is an antiholomorphic morphism
$M|_{U_i} \to M|_{U_{\overline{\i}}}$ (resp. $N|_{U_i} \to
M|_{U_{\overline{\i}}}$, resp. $N|_{U_i} \to N|_{U_{\overline{\i}}}$)
which lifts $c_B$. Since $c_X$ is an involution, we have for every $i
\in J$, $a_{\overline{\i}} \circ a_i = d_{\overline{\i}} \circ d_i
\in {\cal O}_B^*|_{U_i}$ and $a_{\overline{\i}} \circ b_i + b_{\overline{\i}} 
\circ d_i
=0 \in {\cal H}om (N,M)|_{U_i}$. Moreover, for $i,j \in J$
such that $U_i \cap U_j \neq \emptyset$, the gluing conditions are the
following : $a_i = \lambda a_j$, $d_i = \lambda d_j$ and $b_i +
\mu_{\overline{\i} \overline{\j}} \circ d_i = \lambda (a_j \circ
\mu_{ij} + b_j)$ where $\lambda \in {\cal O}_B^*|_{U_i \cap U_j}$.

Now let $Y$ be the complex analytic manifold of dimension $n+1$
defined as the gluing of the charts $\C \times P(M \oplus N)|_{U_i}$, $i \in
J$, with change of charts given by the maps :
\begin{eqnarray*}
\C \times P(M \oplus N)|_{U_i} & \to & \C \times P(M \oplus N)|_{U_j} \\
(t,x,(m:n)) & \mapsto &  (t,x,\left[ \begin{array}{cc}
1&t \mu_{ij}\\
0&1
\end{array} \right]
\left( \begin{array}{c}
m\\
n
\end{array} \right))
=(t,x, (m+t \mu_{ij} n : n)).
\end{eqnarray*}
The projection on the first coordinate defines a holomorphic
submersion $\pi : Y \to \C$. The surface $\pi^{-1} (0)$ is isomorphic
to the ruled manifold $P(M \oplus N)$, whereas, as soon as
$t \in \C^*$, the fiber $Y_t = \pi^{-1} (t)$ is isomorphic to the
ruled manifold $X=P(E)$. Such an isomorphism $\psi_t : Y_t \to X$ is
given in the charts $P(M \oplus N)|_{U_i}$, $i \in J$, by :
\begin{eqnarray*}
P(M \oplus N)|_{U_i} & \to & P(M \oplus N)|_{U_i} \\
(x,(m:n)) & \mapsto &  (x, (m : t n)).
\end{eqnarray*}
Denote by $c_Y$ the real structure on $Y$ defined on charts $\C \times
P(M \oplus N)|_{U_i}$ by :
\begin{eqnarray*}
\C \times P(M \oplus N)|_{U_i} & \to & \C \times P(M \oplus N)|_{U_{\overline{\i}}} \\
(t,x,(m:n)) & \mapsto &  (\overline{t},c_B (x),\left[ \begin{array}{cc}
a_i&\overline{t} b_i\\
0&d_i
\end{array} \right]
\left( \begin{array}{c}
m\\
n
\end{array} \right)).
\end{eqnarray*}
This real structure satisfies $\pi \circ c_Y = conj \circ \pi$ where
$conj$ is the complex conjugation on $\C$. Moreover, when $t \in
\R^*$, $\psi_t$ gives an isomorphism between the real ruled manifolds
$(Y_t , c_Y|_{Y_t})$ and $(X , c_X)$. Hence, the restriction of $\pi :
Y \to \C$ over $\Delta \subset \C$ defines a real deformation and  
the real ruled manifold $P(M \oplus N)$ is thus deformation equivalent to
$(X , c_X)$.
By iteration of this process to $P(N) \subset P(M \oplus N)$, we obtain the
desired result. $\square$ \\

From now on, we will assume that $X=P(L_1 \oplus \dots \oplus L_k \oplus F_1 
\oplus \dots \oplus F_l)$ and that $c_X$ satisfies the conditions of
proposition \ref{propone}. It follows from proposition \ref{propone} that 
this can be done without changing the deformation class of $(X , c_X)$.

\begin{prop}
\label{proptwo}
Let $X=P(L_1 \oplus \dots \oplus L_k \oplus F_1 
\oplus \dots \oplus F_l)$ be a ruled manifold over $B$ and
$c_X$ be a real structure satisfying the conditions of proposition
\ref{propone}. Assume that the real part of $B$ is non-empty, then :

1. If $k \geq 1$, without changing the deformation class of $(X , c_X)$, we
can assume that $l=0$.

2. If $k=0$, without changing the deformation class of $(X , c_X)$, we
can assume that the real ruled surfaces $P(F_j) \subset X$, $j \in \{ 1, 
\dots , l\}$, are obtained from a same real decomposable ruled surface 
after at most one elementary transformation on each
component of its real part if $\R X \neq \emptyset$ and
after at most one couple of elementary transformations on complex conjugated 
points if $\R X = \emptyset$.
\end{prop}

\begin{lemma}
\label{lemmadiv}
Let $L \in \Jac (B)$ be a complex line bundle belonging to the same component
of the real part of $(\Jac (B) , -c_B^*)$ as the trivial line bundle. Then
there exists a divisor $D$ associated to $L$ such 
that $c_B (D) = -D$. The same conclusion holds for any line bundle in the
real part of $(\Jac (B) , -c_B^*)$ when $\R B$ is empty.
\end{lemma}
(In this lemma, the Jacobian $\Jac (B)$ is identified with the degree zero
part of the Picard manifold of $B$.)\\

{\bf Proof :}

Let $L$ be a line bundle belonging to the same component
of the real part of $(\Jac (B) , -c_B^*)$ as the trivial line bundle. Then
there exists a line bundle $\widetilde{L}$ in this component such that
$\widetilde{L} \otimes \widetilde{L} = L$. Let $\widetilde{D}$ be a divisor
associated to $\widetilde{L}$. Then $-c_B (\widetilde{D})$ is also
associated to $\widetilde{L}$ and $D = \widetilde{D} -c_B (\widetilde{D})$ is
suitable. Now assume that $\R B$ is empty. Remember that if the genus
of $B$ is even (resp. is odd), then the real part of $(\Jac (B) , -c_B^*)$ 
is connected (resp. has two connected components), see \cite{GH}, proposition
$3.3$. Assume thus that the genus of $B$ is odd and pick up a point $x$ in
$B$. The line bundle $L(x - c_B (x))$ associated to the divisor $x - c_B (x)$
does not belong to the same component of the real part of $(\Jac (B) , -c_B^*)$
as $L_0$. Indeed the quotients of the real ruled surfaces $(P(L(x - c_B (x))
\oplus L_0),c_X^\pm)$ by the real structures $c_X^\pm$ are not spin
whereas they are for any real ruled surfaces $(P(L
\oplus L_0),c_X^\pm)$ for $L$ in the same component as $L_0$, see
proposition \ref{propempty}. Making the tensor product with the line bundle 
$L(x - c_B (x))$ if necessary, we get the result for any line bundle of
the real part of $(\Jac (B) , -c_B^*)$ in this case. $\square$\\

{\bf Proof of proposition \ref{proptwo} :}

To begin with, let us prove the first part of proposition \ref{proptwo}.
For this, we suppose that the integer $l$ given by proposition \ref{propone}
is non-zero.
Since $k \geq 1$, $X$ admits some real sections and so the projection
$\R X \to \R B$ is surjective. It follows that every real fiber of $X$ has a 
non-empty real part and thus is isomorphic to $\C P^{n-1}$ equipped with 
the standard complex conjugation. Let $j \in \{ 1, \dots , l \}$, the ruled 
surface $P (F_j)$ in $X$ is real and has no holomorphic section. Every real 
fiber
of this surface is a real line in a real fiber of $X$, it thus has real points
and so the real part of $P (F_j)$ has as many components as $\R B$. It
follows then from \cite{Wels}, proposition 3.10, that $P (F_j)$ can be
obtained from a decomposable ruled surface after a finite number of
elementary transformations on real or complex conjugated points. This
decomposable ruled surface is of the form $(P (M \oplus L_0), c_X^\pm)$,
where $M \in \Pic (B)$ and $c_B^* (M) = M^*$ (see \cite{Wels}). Since
the real part of $P (F_j)$ has as many components as $\R B$, it follows from
\cite{Wels}, lemma 3.2, that $M$ belongs to the same component of the
real part of $(\Jac (B) , - c_B^*)$ as $L_0$. From lemma \ref{lemmadiv}, we 
deduce that there exists a divisor $D$ associated to $M$ such that 
$c_B (D) = -D$. So there exists some points $x_i \in B$ and integers $n_i$
such that $D = \sum_i n_i (x_i - c_B (x_i))$. From corollary \ref{corelem}
then follows that $(P (M \oplus L_0), c_X^\pm)$ is obtained from
$(B \times \C P^1 , c_B \times conj)$ after a finite number of elementary
transformations in complex conjugated points. It suffices indeed to make
elementary transformations on the section $P (M)$ over the points
$c_B (x_i)$ and on the section $P (L_0)$ over the points
$x_i$ to pass from $(P (M \oplus L_0), c_X^\pm)$ to $(B \times \C P^1 , c_B 
\times conj)$, since $c_B 
\times conj$ is the only real structure on $B \times \C P^1$ fibered over
$c_B$ and with non-empty real part. 
In conclusion, the real ruled manifold $(X, c_X)$ is obtained from a manifold
$(Y=P(L_1 \oplus \dots \oplus L_k \oplus N \oplus N  \oplus
F_2 \oplus \dots \oplus F_l), c_Y)$ where $c_Y$ satisfies the conditions of
proposition \ref{propone}, after a finite number of elementary transformations
on real and complex conjugated points. Using a real deformation, these points 
can be brought to the real section $P(L_1) \subset Y$. This gives a new real
ruled manifold in the same deformation class as $(X, c_X)$ but with a lower
$l$. After an iteration of this process, we get the result.

Now, let us prove the second part of proposition \ref{proptwo}. 
As in the first part, we can prove that every real ruled
surface $P (F_j)$ is obtained from a decomposable real ruled surface 
$(P (M_j \oplus L_0), c_X^\pm)$ after a finite number of
elementary transformations on real or complex conjugated points. Moreover, 
the line bundle $M_j$ satisfies $c_B^* (M_j) = 
M_j^*$. Also, using
the same trick as in the begining, we deduce that the real parts of
these surfaces project exactly on the same subset of $\R B$. Thus, using the
terminology of \cite{Wels}, we obtain that the partition of $\R B$ in
two elements associated to the bundles $M_j $ are the same. From 
\cite{Wels}, lemma 3.2, it follows that the bundles $M_j$ are in the same 
component
of the real part of $(\Jac (B) , - c_B^*)$ and thus are of the form
$L \otimes \widetilde{M}_j$ with $L, \widetilde{M}_j \in \Jac (B)$,
$c_B^* (L)=L^*$ and $\widetilde{M}_j$ in the same real component of
$(\Jac (B) , - c_B^*)$ as $L_0$. Hence we deduce as before that every 
surface $P (F_j)$ is obtained from a same
decomposable  real ruled surface $(P(L \oplus L_0), c_X^\pm)$ after a finite
number of elementary transformations on real or complex conjugated points.
Now if $\R X$ is non-empty, from corollary \ref{corelem} every couple of  
complex conjugated points can be
deformed into a double real point and then into two real points. Moreover,
 every
couple of real points lying in a same real component of $P(L \oplus L_0)$ can 
be
removed, making the elementary transformation at the first point and bringing
the second on the image of the contracted fiber. This gives the result if 
$\R X$ is non-empty, and otherwise the result is deduced from the
fact that every pair of elementary transformations done on
couples of complex conjugated points of $P(L \oplus L_0)$ can be similarly 
removed. $\square$

\begin{rem}
\label{remproptwo}
From this proof follows that the common decomposable real ruled surface can 
be chosen of the form $(P(L \oplus L_0), c_X^\pm)$ with
$L \neq L^*$ if $L$ is non-trivial.
\end{rem}

We can now prove theorem \ref{theoodd} when the bases of the manifolds
have non-empty 
real parts (the case of empty real part is postponed to next subparagraph).
We separate the cases of odd and even dimensions.\\

{\bf Proof of theorem \ref{theoodd} in odd dimension:}

Let $(X , c_X)$ be a real ruled manifold of odd dimension $n$ and base
$(B , c_B)$ with non-empty real part. Since the integer
$k$ given by proposition \ref{propone} has the same parity as $n$, it is
non-zero. From proposition \ref{proptwo}, we can thus assume that
$k=n$ and $X=P(L_1 \oplus \dots \oplus L_n)$, every section $ P(L_i)$ being
real. We can assume that $L_1$ is trivial, making the tensor product by
$L_1^*$ otherwise. For every $i \in \{ 2 ,\dots , n\}$, the ruled surface
$X_i = P (L_1 \oplus L_i)$ is real. The normal bundle of $ P(L_1)$ in 
$X_i$ is real and isomorphic to $L_i$. There thus
exists a divisor $D_i$ on $B$ such that $c_B (D_i) = D_i$ and $L_i = L (D_i)$.
Hence, from corollary \ref{corelem}, the real ruled manifold $(X , c_X)$ is
 obtained from $(B \times \C P^{n-1} , c_B \times conj)$ after a finite number
of elementary transformations on real or complex conjugated points (remember
that $c_B \times conj$ is the only real structure on $B \times \C P^{n-1}$
fibered over $c_B$, since $n$ is odd). We can assume that on every component
of $\R B \times \R P^{n-1}$, the number of such points is even. Indeed, we
can add $n$ to this number otherwise, making an elementary transformation
on every section $P(L_i)$ over a same point of $\R B$, which does not
affect the isomorphism class of $X$ from corollary \ref{corcont} and 
proposition \ref{proppe}. So from
corollary \ref{corelem} we can assume that all these points are complex 
conjugated, deforming every couple of real points on a same component to a 
double real point and then to two complex conjugated points. Finally, from
corollary \ref{corcont}, the number of elementary transformations done
on complex conjugated points can be reduced modulo $2n$ to an
even number in between zero and $2n - 2$. Indeed, $n$ couple of elementary 
transformations done on complex conjugated points can be removed bringing 
these 
points to the $n$ sections $P(L_i)$ over a same pair of complex conjugated
points of $B$. Note that the ruled manifold obtained from $B \times \C P^{n-1}$
after $2d$ elementary transformations has degree $2d$, so when $d$ ranges
from $0$ to $n-1$, this gives ruled manifolds non deformation equivalent
to each other. 

The proof of theorem \ref{theoodd} is now clear. Let $(X_1 , c_{X_1})$
and $(X_2 , c_{X_2})$ be two real ruled manifolds of odd dimension $n$ with
same degree $d \in \Z/n \Z$ and bases $(B_1 , c_{B_1})$
and $(B_2 , c_{B_2})$ of same topological type. We know that $(X_1 , c_{X_1})$
(resp. $(X_2 , c_{X_2})$) is in the same deformation class as the manifold
obtained from
$(B_1 \times \C P^{n-1} , c_{B_1} \times conj)$ (resp. $(B_2 \times \C P^{n-1}
 , c_{B_2} \times conj)$) after $d$ elementary transformations done on complex
conjugated points. Since $(B_1 , c_{B_1})$
and $(B_2 , c_{B_2})$ have same topological type, from \cite{Nat2} follows
that they are deformation equivalent and there exists a path $(B_t , c_{B_t})$
joining these curves. The path $(X_t , c_{X_t})$, where $(X_t , c_{X_t})$
is obtained from $(B_t \times \C P^{n-1} , c_{B_t} \times conj)$ after
making $d$ elementary transformations on complex
conjugated points, gives a real deformation in between $(X_1 , c_{X_1})$
and $(X_2 , c_{X_2})$. $\square$\\

{\bf Proof of theorem \ref{theoodd} in even dimension :}

Let $(X , c_X)$ be a real ruled manifold of even dimension $n$ and base
$(B , c_B)$ with non-empty real part. Without
changing the deformation class of $(X , c_X)$, we can assume that 
$X=P(L_1 \oplus \dots \oplus L_k \oplus F_1 \oplus \dots \oplus F_l)$ and that
$c_X$ satisfies the conditions of proposition \ref{propone}. From proposition 
\ref{proptwo}, either $k$ or $l$ vanishes.\\

{\bf First case :} $l=0$. In this case, we proceed as in the proof of 
theorem \ref{theoodd} in odd dimension. It follows that the manifold 
$(X , c_X)$ is obtained
from $(B \times \C P^{n-1} , c_B \times conj)$ after a finite number
of elementary transformations done on real or complex conjugated points. 
As in the proof of proposition 
\ref{proptwo},  the number of elementary
transformations can be reduced to zero or one on every component of
$\R B \times \R P^{n-1}$. 
As in the proof of theorem \ref{theoodd} in odd dimension, the number
 of elementary
transformations done on complex conjugated points can be reduced modulo $2n$
to an even number in between $0$ and $2n-2$ and even in between $0$ and $n-2$
since $\R B \neq \emptyset$ and so complex conjugated points can be 
brought to a same real fiber of $X$ and reduced modulo $n$. 
Let $(X_1 , c_{X_1})$
and $(X_2 , c_{X_2})$ be two real ruled manifolds of even dimension $n$ with
same degree $d \in \Z/n \Z$ and same topological type.
Then the number of components of $\R B_1 \times \R P^{n-1}$ (resp. $\R B_2 
\times \R P^{n-1}$) on which is done one elementary transformation is given
by the number of non-orientable component of $\R X_1$ (resp. $\R X_2$), so
it is the same for $X_1$
and $X_2$. Since $(B_1 , c_{B_1})$
and $(B_2 , c_{B_2})$ have same topological type, from \cite{Nat2} follows
that there exists a path $(B_t , c_{B_t})$
joining these curves. Moreover, this path can be chosen so that the component 
of $\R B_1$ over which are the non-orientable components of $\R X_1$ are
mapped to the components of $\R B_2$ over which are the non-orientable 
components of $\R X_2$. This follows from the presentation in \cite{Nat2}
of a real algebraic curve as the gluing of a Riemann surface with
boundary with its conjugate, the gluing maps being either identity or
antipodal. Hence, we can assume that $(X_1 , c_{X_1})$
and $(X_2 , c_{X_2})$ are obtained from a same real ruled manifold $(Y, c_Y)$
after making an even number of elementary
transformations on  complex conjugated points, these numbers being less
than $n-1$. Since by hypothesis these manifolds have same degree, these numbers
are the same for $(X_1 , c_{X_1})$
and $(X_2 , c_{X_2})$ and these two manifolds are deformation equivalent.\\

{\bf Second case :} $k=0$. Then, from proposition \ref{proptwo}, 
without changing the deformation class of $(X, c_X)$ and making a finite
number of elementary transformations on real points if necessary, we can 
assume that there exists $L \in \Pic (B)$ such that $c_B^* (L) = L^*$ and
for $1 \leq j \leq \frac{1}{2} n$, $P(F_j) = P (L \oplus L_0)$. Denote by
$M_j$ the line bundle such that $F_j = (M_j \otimes L) \oplus M_j$. We can
assume that $M_1$ is trivial and we denote by $c_1$ the real structure on
$P(F_1)$ induced by $c_X$. The real structure $c_X$ induces a real structure
on the normal bundle $N_j$ of $P(F_1)$ in $P(F_1 \oplus F_j)$. Thus
$c_1^* (N_j) = N_j$. But from lemma \ref{lemmanorm}, we know that
$$N_j = p^* (F_j) \otimes {\cal O} (D_L) = p^* ((M_j \otimes L) \oplus M_j)
\otimes {\cal O} (D_L) = ({\cal O} (D_L) \oplus {\cal O} (D_{L_0}))
\otimes p^* (M_j).$$ 
So $c_1^* (N_j) = N_j$ implies that $c_B^* (M_j) = M_j$. Since $\R B \neq 
\emptyset$, this implies that there
exists a divisor $D_j$ on $B$ associated to $M_j$ such that $c_B (D_j) = D_j$.
We then deduce from corollary \ref{corelem} and proposition \ref{propconj}
that $(X, c_X)$ is obtained from $(Y=P ((L \oplus L_0)^\frac{n}{2}) , c_Y^\pm)$
after a finite number of elementary transformations on real or complex 
conjugated points (we can indeed assume that $L \neq L^*$ if $L$ is 
non-trivial, see remark \ref{remproptwo}). As in the first case, the number
 of elementary
transformations done on real points can be reduced to $0$ or $1$ for each 
component of $\R Y$ and on complex conjugated points, they can be reduced 
modulo $n$ to an even number in between $0$ and $n-2$.
We conclude exactly as in the first case. 
$\square$

\subsection{When the base has empty real part}
\label{parempty}

The aim of this paragraph is to prove theorem \ref{theoodd} assuming
the bases of the manifolds have empty real parts. This is the only remaining
case to consider, after \S \ref{parnonempty}.

\begin{lemma}
\label{lemmasect}
Let $(X= P(L_0 \oplus L \oplus L ), c_X )$ be a real ruled manifold having a 
base $(B , c_B)$ with empty real part. Assume that the real structure $c_X$
fixes the section $P(L_0 ) \subset X$ and exchanges the two sections
$P(L ) \subset X$. Then $X$ is obtained from $(B \times \C P^3 , c_B 
\times \conj)$ after a finite number of elementary transformations done
on real or complex conjugated points.
\end{lemma}

{\bf Proof :}

Denote by $D_L$ the divisor $P(L) \subset P(L \oplus L)$, where $P(L)$ is one
of the two sections of $P(L \oplus L)$, and by $c_1$ the real structure
of $P(L \oplus L)$ induced by $c_X$. The normal bundle of $P(L \oplus L)$
in $X$ is real and isomorphic to $p^* (L^*) \otimes {\cal O} (D_L)$ from
lemma \ref{lemmanorm}. Thus, $c_1^* (p^* (L^*) \otimes {\cal O} (D_L)) =
p^* (L^*) \otimes {\cal O} (D_L)$. Now $c_1^* (p^* (L^*) \otimes {\cal O} 
(D_L)) = p^* (c_B^* (L^*)) \otimes {\cal O} (D_L)$, so that $c_B^* (L) = L$.
Let $D$ be a divisor associated to $L$. Then $c_B (D)$ is also associated
to $L$ and one can write $(X= P(L_0 \oplus L(D) \oplus L(c_B (D)) )$. Making
elementary transformations on $P(L_0)$ over the points of $D + c_B (D)$,
on $P(L(D))$ over the points of $c_B (D)$ and on $P(L(c_B (D))$ over the
points of $D$, we obtain the real ruled manifold 
$(B \times \C P^3 , c_B \times \conj)$. Hence the result. $\square$

\begin{prop}
\label{prop3}
Let $X=P(L_1 \oplus \dots \oplus L_k \oplus F_1 
\oplus \dots \oplus F_l)$ be a ruled manifold over $B$ and
$c_X$ be a real structure satisfying the conditions of proposition
\ref{propone}. Assume that the real part of $B$ is empty, then :

1. If $k \geq 1$, without changing the deformation class of $(X , c_X)$, we
can assume that $l=0$.

2. If $k=0$, without changing the deformation class of $(X , c_X)$, we
can assume that the real ruled surfaces $P(F_j) \subset X$, $j \in \{ 1, 
\dots , l\}$, are obtained from the trivial ruled surface 
after at most one couple of elementary transformations done on complex 
conjugated points. 
\end{prop} 

{\bf Proof :}

Thanks to lemma \ref{lemmadiv}, we can proceed as in the proof of proposition
\ref{proptwo} to get that all the real ruled surfaces $P(F_j)$ are 
obtained from the trivial ruled surface after some elementary
transformations done on complex conjugated points. If $k \geq 1$, all the
elementary transformations can be brought to the section $P(L_1)$ without
changing the real deformation class of $X$. If $k=0$, every two such couples
can be cancelled, making the first two elementary transformations and bringing
the two remaining ones on the images of the contracted fibers.
This already proves the second part of proposition \ref{prop3}. Now if
$k,l \geq 1$, we can assume that $L_1$ is the trivial line bundle $L_0$. From
what has just be done, we know that $F_j = M_j \oplus M_j$ for some
line bundle $M_j$ over $B$. From lemma \ref{lemmasect}, we deduce that the
real ruled surfaces $P(F_j) = P(M_j \oplus M_j)$ have a real section, which 
contradicts the hypothesis (see proposition \ref{propone}). 
$\square$

\begin{prop}
\label{propempty}
Let $n \geq 2$ be an even integer and $(B , c_B)$ be a real compact
irreducible algebraic curve with empty real part. Let $c_0$ be the real 
structure with empty real part on $\C P^{n-1}$.

1) If $g(B)$ is odd, then the real ruled manifolds $(B \times \C P^{n-1} ,
c_B \times \conj)$ and  $(B \times \C P^{n-1} , c_B \times c_0)$ are
deformation equivalent.

2)Let $(X , c_X)$ be the real ruled manifold obtained from $(B \times 
\C P^{n-1} , c_B \times \conj)$ after a couple of elementary transformations
done along two projective subspaces of dimension $\frac{1}{2} n -1$ 
belonging to
two complex conjugated fibers. Then the quotients $X/c_X$ and $B \times 
\C P^{n-1} / c_B \times \conj$ are not homeomorphic to each other.
\end{prop}

\begin{cor}
\label{corempty}
Let $n \geq 2$ be an even integer and $(B , c_B)$ be a real compact
irreducible algebraic curve with empty real part and even genus. Then
the quotients $B \times \C P^{n-1} / c_B \times \conj$ and $B \times 
\C P^{n-1} / c_B \times c_0$ are not homeomorphic to each other.
\end{cor}

{\bf Proof :}

From proposition \ref{propempty}, we know that the real ruled manifold 
$(X , c_X)$ obtained from $(B \times 
\C P^{n-1} , c_B \times \conj)$ after a couple of elementary transformations
done along two projective subspaces of dimension $\frac{1}{2} n -1$ of
two complex conjugated fibers has a quotient non homeomorphic to
$B \times \C P^{n-1} / c_B \times \conj$. Let us write $B \times \C P^{n-1} =
P(L_0^\frac{n}{2} \oplus L_0^\frac{n}{2})$, where the two sub-ruled manifolds
$P(L_0^\frac{n}{2})$ are exchanged by $c_B \times \conj$. We can perform
the elementary transformations along these subspaces over two complex
conjugated points $x$ and $c_B (x)$ of $B$. From corollary \ref{corelem} we
then deduce that $X = P((L(x - c_B(x)) \oplus L_0)^\frac{n}{2})$ and from
proposition \ref{propconj}, we deduce that $c_X$ is one of 
the two real structures $c_X^\pm$. Now since the genus of $B$ is even, the real
part of $(\Jac (B) , -c_B^*)$ is connected and hence $(X , c_X)$ deforms
onto $B \times \C P^{n-1}$ equipped with a real structure having a
quotient homeomorphic to $X / c_X$. From proposition \ref{propempty}
it cannot be $c_B \times \conj$, hence the result. $\square$ \\

{\bf Proof of proposition \ref{propempty} :}

The proof of the first part of this proposition is analogous to the one of
proposition $2.8$ of \cite{Wels}, so we will give only a sketch of it.
Without changing the deformation class of the manifolds, we can assume that the
base $(B, c_B)$ is the curve constructed in corollary $2.8$ of \cite{Wels}.
Let $L$ be the line bundle given by this corollary, it satisfies $c_B^* (L) =
L=L^*$. Moreover, there exist a divisor $D$ associated to $L$, an automorphism
$\varphi$ of $B$, a meromorphic function $f_D$ on $B$ such that
$\overline{f \circ c_B} = f$, $\Div (f) = D + c_B (D)$ and 
$f \circ \varphi =  f$, as well as a meromorphic function $g$ on $B$ such 
that $\Div (g) = \varphi (D) - D$ and $g (\overline{g \circ c_B}) = -1$
(see \cite{Wels}).
We define, as in the proof of proposition $2.6$ of \cite{Wels}, an
automorphism $\Phi$ of $P((L \oplus L_0)^\frac{n}{2})$ defined over the
open set $U_0 = B \setminus D$, by :
$$\begin{array}{rcl}
U_0  \times \C P^{n-1} & \to & U_0 \times \C P^{n-1}\\
(x,(z_1^j : z_0^j)) & \mapsto & 
(\varphi (x), (g \circ \varphi (x) z_1^j : z_0^j)), 
\end{array}$$
where $j \in \{ 1, \dots , \frac{n}{2} \}$. This automorphism conjugates the
two real structures $c_X^+$ and $c_X^-$ on $X$. Now if $L$ is in the same
component of the real part of $(\Jac (B) , -c_B^*)$ as the trivial line
bundle, the real ruled manifolds $(B \times \C P^{n-1} ,
c_B \times \conj)$ and  $(B \times \C P^{n-1} , c_B \times c_0)$ are in
the same deformation class as $(P((L \oplus L_0)^\frac{n}{2}) , c_X^\pm)$
and we are done. Otherwise, we perform an elementary transformation along
the subspace $P(L^\frac{n}{2}) \subset P((L \oplus L_0)^\frac{n}{2})$ over
a point $x$ of $B$, and an elementary transformation along
the subspace $c_X^+ (P(L^\frac{n}{2})) = c_X^- (P(L^\frac{n}{2}))$ over the
point $c_B (x)$. The ruled manifold obtained is of the form
$P((M \oplus L_0)^\frac{n}{2})$ with $M$ in the same
component of the real part of $(\Jac (B) , -c_B^*)$ as $L_0$. The real 
structures $c_X^\pm$ lift to the real structures $c_Y^\pm$, and the 
automorphism $\Phi$ lifts to an automorphism of $Y$ which conjugates them
(see lemma $3.14$ of \cite{Wels}). Hence we conclude as before.

Now let us prove the second part of the proposition. We will prove that there
is no $\Z /2\Z$-equivariant diffeomorphism in between
$(B \times \C P^{n-1} , c_B \times \conj)$ and $(X , c_X)$. Let $H$ be
a real hyperplane of $(\C P^{n-1} , \conj)$ and $s = B \times H$. Then $s$ is
a real divisor of $(B \times \C P^{n-1} , c_B \times \conj)$ and $[s]^n = 0
\in H_0 (B \times \C P^{n-1} ; \Z)$. If $s'$ is any other real 
$(2n-2)$-cycle of $(B \times \C P^{n-1} , c_B \times \conj)$, then $[s'] =
[s] + 2k[f]$ where $[f]$ is the homology class of a fiber. Hence
$[s']^n = [s]^n +2kn [f][s]^{n-1} = \pm 2kn = 0 \mod (2n)$. Now, without
changing the diffeomorphism class of $X/c_X$, we can assume that the blown up
projective subspaces of dimension $\frac{1}{2} n - 1$ are real once projected
onto the second factor $(\C P^{n-1} , \conj)$, so that $X = 
P(L (x + c_B (x))^\frac{n}{2} \oplus L_0^\frac{n}{2})$ and $c_X$ is the 
standard real structure on this manifold. Denote by $s_X$ the real divisor 
on $X$
associated to the dual of the tautological line bundle ${\cal O} (-1)
\subset p^* (L (x + c_B (x))^\frac{n}{2} \oplus L_0^\frac{n}{2})$. From
\cite{Schm}, p. $215$, we see that $[s_X]^n=- \deg 
(L (x + c_B (x))^\frac{n}{2} \oplus L_0^\frac{n}{2}) = -n \neq 0 \mod(2n)$,
hence the result. $\square$ \\

{\bf Proof of theorem \ref{theoodd} when the base has empty real part :}

Let $(X,c_X)$ be a real ruled manifold of dimension $n$ and base $(B,c_B)$
with empty real part. If $n$ is odd, from propositions \ref{propone} and
\ref{prop3}, we know that without changing the deformation class of $(X,c_X)$,
we can 
assume that $X=P(L_1 \oplus \dots \oplus L_n)$, all the sections
$P(L_j)$ being real. We can also assume that $L_1$ is the trivial line bundle.
Then, for $2 \leq j \leq n$, the normal bundle of $P(L_0)$ in $P(L_0 \oplus
L_j)$ is real and isomorphic to $L_j$. Thus there exists a real divisor $D_j$
associated to $L_j$ and from corollary \ref{corelem}, $(X,c_X)$ is
obtained from $(B \times \C P^{n-1}, c_B \times \conj)$ after a finite number 
of elementary transformations done on complex conjugated points. Note
that since $n$ is odd, $c_B \times \conj$ is the only real structure on
$B \times \C P^{n-1}$ fibered over $c_B$. Now from corollary \ref{corelem},
without changing the deformation class of $(X,c_X)$, $n$ couples of elementary
transformations done on complex conjugated points can be removed, bringing
these points on the $n$ sections over a same couple $(x,c_B (x))$ of $B
\times B$. Hence we can assume that $(X,c_X)$ is
obtained from $(B \times \C P^{n-1}, c_B \times \conj)$ after an even
number of elementary transformations in between $0$ and $2n-2$. Since
$n$ is odd, this number is prescribed by the degree of $X$, and since
there is only one deformation class of smooth real irreducible compact
algebraic curve with empty real part (see \cite{Nat2}), the result is proved
in odd dimension.

Now if $n$ is even and the integer $k$ given by proposition \ref{prop3} is
non-zero, the same proof shows that without changing the deformation class 
of $(X,c_X)$, we can assume that it is obtained from $B \times \C P^{n-1}$
equipped with one of its real structures after an even number of elementary 
transformations done on complex conjugated points, in between $0$ and $2n-2$.
In particular, the degree of $X$ is even.
If $g(B)$ is even, then from corollary \ref{corempty}, the two real 
structures on $B \times \C P^{n-1}$ have non-diffeomorphic quotients, and from
proposition \ref{propempty}, one passes from the deformation class of one
of these real structures to the deformation class of the other one making
$n$ elementary transformations on complex conjugated points. Thus we can assume
that the number of elementary transformations necessary to obtain $(X,c_X)$
from $B \times \C P^{n-1}$ with one of its real structures is even in between
$0$ and $n-2$. This number is then prescribed by the degree of $X$, and the
real structure on $B \times \C P^{n-1}$ by the topology of the quotient
$X / c_X$. If $g(B)$ is odd, then from proposition \ref{propempty}, the
two real structures on $B \times \C P^{n-1}$ are in the same deformation
class, hence the  number modulo $n$ of elementary transformations necessary 
to obtain $(X,c_X)$ from $B \times \C P^{n-1}$ is prescribed by the degree of 
$X$, and then the total number by the topology of the quotient $X / c_X$. We
can then conclude as before.

It thus only remains to consider the case when $n$ is even, but the integer
$k$ given by proposition \ref{prop3} vanishes. From proposition \ref{prop3},
we can assume that $(X,c_X)$ is obtained from a real ruled manifold
$P(F_1 \oplus \dots \oplus F_l)$ after some elementary transformations 
done on complex conjugated points. Moreover, all the ruled surfaces $P(F_j)$
are real, and every bundle $F_j$ is of the form $M_j \oplus M_j$ for
some $M_j \in \Pic (B)$. We can then assume that $M_1$ is trivial and since
the normal bundle of $P(F_1)$ in $P(F_1 \oplus F_j)$ is real, we deduce
from lemma \ref{lemmanorm} that $c_B^* (M_j) = M_j$, as in the proof of
lemma \ref{lemmasect}. Let $D_j$ be a divisor associated to $M_j$. Then
$c_B (D_j)$ is also associated to $M_j$ and we can write
$P(F_1 \oplus \dots \oplus F_l) = P(L_0 \oplus L_0 \oplus L(D_2)
 \oplus L(c_B (D_2)) \oplus \dots \oplus L(D_l) \oplus L(c_B (D_l)))$.
Making
elementary transformations on $P(L_0 \oplus L_0)$ over the points of
 $D_2 + c_B (D_2) + \dots + D_l + c_B (D_l)$,
on $P(L(D_j))$ over the points of $c_B (D_j)$ and on $P(L(c_B (D_j))$ over the
points of $D_j$ (with appropriate multiplicities), we obtain the real ruled 
manifold 
$B \times \C P^{n-1}$ with one of its real structures. Hence once more,
$(X,c_X)$ is in the same deformation class as a manifold obtained from
$B \times \C P^{n-1}$ with one of its real structures after a finite number
of elementary transformations done on complex conjugated points and we
conclude as before. $\square$

\addcontentsline{toc}{part}{\hspace*{\indentation}Bibliographie}


\bibliography{ruledmfd}
\bibliographystyle{abbrv}

\noindent Ecole Normale Sup\'erieure de Lyon\\
Unit\'e de Math\'ematiques pures et appliqu\'ees\\
$46$, all\'ee d'Italie\\
$69364$, Lyon C\'edex $07$\\
(FRANCE)\\
e-mail : {\tt jwelschi@umpa.ens-lyon.fr}

\end{document}